\documentclass[11pt]{article}
\usepackage{amsmath}
\usepackage{amssymb}
\usepackage{latexsym}
\usepackage{amsthm}

\usepackage{graphicx}
\usepackage{amsmath,amsfonts,amssymb}

\usepackage{url}
\usepackage{amsmath}
\usepackage{graphicx}
\usepackage{multirow}
\usepackage{wrapfig}
\usepackage{citesort}
\usepackage{epsfig}
\usepackage{color}
\def\Tilde{\widetilde}

\def\ra{\rangle}
\def\la{\langle}
\def\ve{\varepsilon}
\def\B{I\!\!B}
\def\h{\hfill\Box}
\def\R{I\!\!R}
\def\N{I\!\!N}
\def\ox{\bar{x}}

\def\h{\hfill\triangle}

\def\O{\Omega}

\def \N{I\!\!N}

\renewcommand{\theequation}{\thesection.\arabic{equation}}
\newcounter{lk}

\begin{document}
\begin{center}
{\bf THE SMALLEST ENCLOSING BALL PROBLEM AND \\THE SMALLEST INTERSECTING BALL PROBLEM:\\ EXISTENCE AND UNIQUENESS OF SOLUTIONS}\\[2ex]
Boris S. Mordukhovich\footnote{Department of Mathematics, Wayne State University,
Detroit, MI 48202, USA (email: boris@math.wayne.edu). The research of these
authors was partially supported by the US National Science Foundation under grant DMS-1007132, by the Australian Research Council under grant
DP-12092508, and by the Portuguese Foundation of Science and Technologies under grant MAT/11109.}, Nguyen Mau Nam\footnote{Department of Mathematics, University of
Texas-Pan American, TX 78539, USA (email: nguyenmn@utpa.edu). The
research of this author was partially supported by the Simons
Foundation under grant \#208785.}, Cristina
Villalobos\footnote{Department of Mathematics, The University of
Texas--Pan American, Edinburg, TX 78539--2999, USA (email:
mcvilla@utpa.edu).}
\end{center}
\small{\bf Abstract:}  In this paper we study the following
problems: given a finite number of nonempty closed subsets of a
normed space, find a ball with the smallest radius that encloses all
of the sets, and find a ball with the smallest radius that
intersects all of the sets. These problems can be viewed as
generalized versions of the \emph{smallest enclosing circle problem}
introduced in the 19th century by Sylvester \cite{syl} which asks
for the circle of smallest radius enclosing a given set of finite
points in the plane. We will focus on the sufficient conditions for
the existence and uniqueness of an optimal solution for each
problem, while the study of optimality conditions and numerical
implementation will be addressed in our next projects.

\medskip
\vspace*{0,05in} \noindent {\bf Key words.} smallest enclosing ball
problem, smallest intersecting ball problem


\newtheorem{Theorem}{Theorem}[section]
\newtheorem{Proposition}[Theorem]{Proposition}
\newtheorem{Remark}[Theorem]{Remark}
\newtheorem{Lemma}[Theorem]{Lemma}
\newtheorem{Corollary}[Theorem]{Corollary}
\newtheorem{Definition}[Theorem]{Definition}
\newtheorem{Example}[Theorem]{Example}
\renewcommand{\theequation}{\thesection.\arabic{equation}}
\normalsize

\section{Introduction and Problem Formulation}
\setcounter{equation}{0}

Let $X$ be a normed space, and let $F \subset X$ be a closed,
bounded, convex set which contains the origin as an interior point.
Given $x\in X$ and $r\geq 0$, we define
\begin{equation*}\label{ball}
D_F(x; r)=x +r F
\end{equation*}
to be  a closed bounded convex set centered about the point $x$ with
radius $r$.

The first problem we propose and study in this paper is called the
\emph{smallest enclosing ball problem} and is stated as follows:
given a nonempty closed constraint set $\O \subset X $ and a
finite number of nonempty closed  bounded subsets $\O_i \subset X$
for $i=1,\ldots,n$, find a point $\ox\in \O$ and the smallest radius
$r\geq 0$ such that the set $D_F(\ox;r)$ contains all of the sets,
that is, $\O_i \subset D_F(\ox;r)$ for $i=1,\ldots,n$.

The second problem under consideration is called the \emph{smallest
intersecting ball problem} and is stated as follows: given a
nonempty closed constraint set $\O \subset X$ and a finite number
of nonempty closed subsets $\O_i \subset X$ for $i=1,\ldots,n$,
find a point $\ox\in \O$ with smallest radius $r\geq 0$ such that
the set $D_F(\ox;r)$ intersects all of the sets $\O_i$.
\vspace*{0.05in}

When  $F$ is the  closed unit ball in the Euclidean plane and the
\emph{target sets} $\O_i$, $i=1,\ldots, n$, are singletons and the
constraint set $\O$ is the whole plane, then both problems reduce to
the \emph{classical smallest enclosing circle problem} introduced by
the English mathematician James Joseph Sylvester (1814--1897) which
asks for the smallest circle that covers a finite number of points
on the plane. After more than a century, the smallest enclosing
circle problem remains active; see \cite{chm,frank,ljc,wel} and the
references therein. The reader is referred to our recent paper
\cite{naj} for a comprehensive study of the convex version of the
smallest intersecting ball problem. The results presented in this
paper and its continuation further our idea of using
variational/nonsmooth analysis and optimization to shed new light on
classical geometry problems.\vspace*{0.05in}

Given a nonempty closed bounded set $Q \subset X$, we define the
\emph{maximal time function}  of $x \in X$ given $Q$ generated by
$F$ as follows
\begin{equation}\label{maxima time function}
C_F(x; Q)=\inf \{ t\geq 0: Q\subset x+tF \}.
\end{equation}
When $F$ is the closed unit ball of $X$, the maximal time function
(\ref{maxima time function}) reduces to the corresponding
\emph{farthest distance function}

\begin{equation*}
M(x; Q)=\sup\{ ||x-\omega|| : \omega\in Q\}.
\end{equation*}
General and generalized differentiation properties of farthest
distance functions can be found, for instance, in \cite{DZ,Lau,WS}.

The \emph{minimal time function} counterpart is  defined below as
\begin{equation}\label{minimal time function}
T_F(x; Q)=\inf \{ t\geq 0: (x+tF)\cap Q\neq \emptyset \},
\end{equation}
where $Q$ needs not necessarily be bounded.  The minimal time
function (\ref{minimal time function}) is more well-known in the
literature; see, e.g. \cite{bmn10} and the references therein. It
becomes the familiar distance function
\begin{equation*}
d(x; Q)=\inf\{ ||x-q||: q\in Q\}
\end{equation*}
when $F$ is the closed unit ball of $X$.\vspace*{0.05in}

In this paper, we will show that under natural assumptions, the
smallest enclosing ball problem can be modeled in terms of an
optimization problem as follows:

\begin{equation}\label{SEC}
\mbox{minimize } C(x) \mbox{ subject to }x\in \O,
\end{equation}
where
\begin{equation}\label{max1}
C(x)=\max\{ C_F(x;\O_i): i=1,\ldots,n\}.
\end{equation}
Similarly, the smallest intersecting ball problem can also be
converted to the following optimization problem:
\begin{equation}\label{SIB}
\mbox{minimize } T(x) \mbox{ subject to }x\in \O,
\end{equation}
where
\begin{equation}\label{max}
T(x)=\max\{ T_F(x;\O_i): i=1,\ldots,n\}.
\end{equation}
The unconstrained versions of these problems are obtained when
$\O=X$.\vspace*{0.05in}

Our goal in this paper and its continuation is to initiate
comprehensive studies of the smallest enclosing ball problem and the
smallest intersecting ball problem using modern tools of variational
analysis and optimization. The main focus of the paper is on
sufficient conditions that guarantee the existence and uniqueness of
an optimal solution for each problem. In Section 2, we provide
important properties of the maximal time function (\ref{maxima time
function}) and then pay  attention to optimality conditions of the
smallest enclosing ball problem. Section 3 is devoted to the
smallest intersecting ball problem counterpart. Along the way, we
point out major differences between these two problems and provide
examples to support the need for the assumptions.  For instance, in
finite dimensional Euclidean space, the smallest enclosing ball
problem usually has a unique optimal solution even if the target
sets are nonconvex, while strict convexity assumptions must be made
to guarantee a unique solution for the smallest intersecting ball
problem.

\section{The Smallest Enclosing Ball Problem}
\setcounter{equation}{0} In this section we initially describe some
properties of the maximal time function (\ref{maxima time function})
and then prove existence and uniqueness of the solution to the
smallest enclosing ball problem in Theorems \ref{existence1} and
\ref{uniqueness1}, respectively.  Finally, we provide some examples
that illustrate the need for the assumptions to guarantee uniqueness
of the solution.

Throughout this section we make the following standing assumptions
unless otherwise noted:

\emph{$X$ is a normed space; $F \subset X $ is a closed, bounded,
convex set that contains the origin as an interior point; the target
sets $\Omega_i$, $i=1,\ldots,n$, are nonempty closed bounded subsets of
$X$; and the constrained set $\Omega$ is a nonempty closed subset of
$X$.}
\\ \\
Let us start with some important properties of the maximal time
function (\ref{maxima time function}). Recall that the Minkowski
function generated by $F$ is given by
\begin{equation}\label{minkowski}
\rho_F(x)=\inf\{ t\geq 0: x\in tF\}.
\end{equation}

The following proposition allows us to represent the maximal time
function (\ref{maxima time function}) in terms of the Minkowski
function (\ref{minkowski}).

\begin{Proposition}\label{p1} Suppose that $Q$ is a nonempty bounded set of $X$. Then the maximal time function (\ref{maxima time function}) has the following representation:
\begin{equation*}
C_F(x; Q)=\sup \{\rho_F(\omega-x): \omega\in Q\}.
\end{equation*}
Moreover, if $F$ is the closed unit ball of $X$, then
\begin{equation*}
C_F(x; Q)=\sup \{||x-\omega||:\omega\in Q\}.
\end{equation*}
\end{Proposition}
{\bf Proof: } Define
\begin{equation*}
\Tilde{C}_F(x;Q)=\sup \{\rho_F(\omega-x): \omega\in Q\}.
\end{equation*}
Fix any $t\geq 0$ such that $Q\subset x+tF$. Then for every
$\omega\in Q$ one has
\begin{equation*}
\omega-x\in tF.
\end{equation*}
Thus $\rho_F(\omega-x)\leq t$. It follows that
\begin{equation*}
\Tilde{C}_F(x;Q)\leq C_F(x; Q).
\end{equation*}
Given any $\ve>0$, one has
\begin{equation*}
\sup \{\rho_F(\omega-x)\; |\; \omega\in Q\}< \Tilde{C}_F(x;Q) +\ve.
\end{equation*}
Thus
\begin{equation*}
\rho_F(\omega-x)<\Tilde{C}_F(x;Q) +\ve \mbox{ for every }\omega\in
Q.
\end{equation*}
By the definition of the Minkowski function, there exists $t\geq 0$
with $t<\Tilde{C}_F(x;Q) +\ve$ and
\begin{equation*}
\omega-x\in tF.
\end{equation*}
Since $F$ is convex and $0\in F$, one has
\begin{equation*}
\omega-x\in (\Tilde{C}_F(x;Q) +\ve)F.
\end{equation*}
This implies
\begin{equation*}
Q\subset x+ (\Tilde{C}_F(x;Q) +\ve)F,
\end{equation*}
and hence $C_F(x;Q)\leq \Tilde{C}_F(x;Q) +\ve$. We finally have
\begin{equation*}
C_F(x;Q)\leq \Tilde{C}_F(x;Q)
\end{equation*}
because $\ve$ is arbitrarily chosen. Thus $C_F(x;Q) =
\Tilde{C}_F(x;Q)$, and the first representation has been proven. In
the case where $F$ is the closed unit ball of $X$, one has
$\rho_F(x)=||x||$. Therefore, the second representation becomes
straightforward. $\h$\vspace*{0.05in}

For a point $x\in X$, the \emph{farthest projection} from $x$ to a
nonempty, closed, bounded set $Q$ with respect to $F$ is defined by
\begin{equation*}
\mathcal{P}_F(x; Q)=\{\omega\in Q: \rho_F(\omega-x)=C_F(x;Q)\}.
\end{equation*}
This set is obviously nonempty when $Q$ is a compact subset of $X$.
The following proposition provides a geometric way to realize the
set.
\begin{Proposition}\label{projection} Let $Q$ be a nonempty, closed, bounded, subset of $X$. Then
\begin{equation*}
\mathcal{P}_F(x; Q)= Q \cap (x+C_F(x; Q)\mbox{bd}(F)),
\end{equation*}
where $\mbox{bd}(F)$ stands for the boundary of $F$.
\end{Proposition}
{\bf Proof: } Fix any $\omega\in \mathcal{P}_F(x; Q)$. Then
$\omega\in Q$ and
\begin{equation*}
\rho_F(\omega-x)=C_F(x; Q).
\end{equation*}
This implies $\omega-x\in C_F(x; Q)F$. When $C_F(x;Q)=0$, it is
obvious that $$\omega\in Q \cap (x+C_F(x; Q)\mbox{bd}(F)) = \{x\}.$$
Suppose $C_F(x;Q)>0$. Then
\begin{equation*}
\rho_F \left(\dfrac{\omega-x}{C_F(x;Q)}\right)=1.
\end{equation*}
 By the well-known property of the Minkowski function, this equality implies that $\dfrac{\omega-x}{C_F(x;Q)}\in \mbox{bd}(F)$, and hence $\omega\in x+C_F(x; Q)\mbox{bd}(F)$. We have shown that
\begin{equation*}
\mathcal{P}_F(x; Q)\subset Q \cap (x+C_F(x; Q)\mbox{bd}(F)).
\end{equation*}
The oppositive inclusion can also be proved similarly.
$\h$\vspace*{0.05in}

Recall that a function $\phi: X\to \R$ is convex on a convex set
$\Omega$ if for every $x,y\in \Omega$ and $t\in (0,1)$, one has
\begin{equation*}
\phi(tx+(1-t)y)\leq t\phi(x)+ (1-t)\phi(y).
\end{equation*}
If this inequality becomes strict for every $x,y\in\Omega$ with $x\neq
y$ and for every $t\in (0,1)$, the function $\phi$ is called
\emph{strictly convex}.

It is clear that if a function $\phi$ is strictly convex on a convex
set $\Omega$, then the problem
\begin{equation*}
\mbox{ minimize } \phi(x) \mbox{ subject to } x\in \Omega
\end{equation*}
cannot have more than one solution.

We also say that a set $\Omega$ is strictly convex if for any $x, y \in
\Omega$ with $x\neq y$ and for any $t\in (0,1)$ we have
\begin{equation*}
tx +(1-t)y\in \mbox{ int }\Omega.
\end{equation*}

\begin{Proposition}\label{convexLipschitz} Let $Q$ be a nonempty  bounded subset of $X$. Then the maximal time function (\ref{maxima time function}) is a finite convex Lipschitz function.
\end{Proposition}
{\bf Proof: } For each $\omega\in Q$, the function
\begin{equation*}
g_\omega(x)=\rho_F(\omega-x)
\end{equation*}
is a convex Lipschitz function since the Minkowski function
$\rho_F(x)$ given in (\ref{minkowski})  is always Lipschitz
continuous with Lipschitz constant $\ell$. Since $\rho_F(0)=0$, one
has $\rho_F(x)\leq \ell ||x||$ for all $x\in X$. It follows from
Proposition (\ref{p1}) that the function
\begin{equation*}
C_F(x;\Omega)=\sup_{\omega\in Q}g_\omega(x)
\end{equation*}
is also a finite convex Lipschitz function under the boundedness
assumption imposed on $Q$. $\h$

\begin{Theorem}\label{existence1} Suppose one of the following holds:\\[1ex]
{\bf (i)} The constraint $\Omega$ is a nonempty compact set.\\[1ex]
{\bf (ii)} $X$ is a reflexive Banach space and the constraint set $\Omega$ is weakly closed.\\[1ex]
Then the smallest enclosing ball problem has a solution. That means
there exist $\ox\in \Omega$ and  $r\geq 0$ such that
\begin{equation*}
\Omega_i\subset \ox+ rF \mbox{ for }i=1,\ldots, n,
\end{equation*}
and for any $x\in \Omega$ and $t\geq 0$ such that $\Omega_i\subset x+tF$ for
$i=1,\ldots,n$, one has $r\leq t$.
\end{Theorem}
{\bf Proof: } Define
\begin{equation*}
\mathcal{E}=\{ t\geq 0: \mbox{there exists } x\in \Omega \mbox{ with }
\Omega_i\subset x+tF \mbox{ for all }i=1,\ldots,n\}.
\end{equation*}
Clearly, $\mathcal{E}\neq \emptyset$ and  is bounded below. Indeed,
fix $\ox\in\Omega$. Since $F$ is convex, $0\in\mbox{int} F$, and $\Omega_i$
is bounded for every $i=1,\ldots,n$, there exists $\bar t>0$ such
that
\begin{equation*}
\Omega_i\subset \ox+\bar tF \mbox{ for all }i=1,\ldots,n.
\end{equation*}
Then $\bar t\in \mathcal{E} \neq \emptyset$. Define $r=\inf
\mathcal{E}$. Let $(t_k)$ be a sequence in $\mathcal{E}$ that
converges to $r$. Let $(x_k)$ be a sequence of $\Omega$ such that
\begin{equation*}
\Omega_i\subset x_k +t_k F \mbox{ for all } i=1,\ldots, n.
\end{equation*}

Now we will consider each case given in the assumptions. In the case
(i) where $\Omega$ is compact, the sequence $(x_k)$ has a subsequence
(without relabeling) that converges to some $\ox\in \Omega$ since $\Omega$
is closed and in this case
\begin{equation}\label{cover}
\Omega_i\subset \ox + r F \mbox{ for all } i=1,\ldots, n.
\end{equation}
Let us now consider case (ii) where $X$ is a reflexive Banach space
and $\Omega$ is weakly closed. It is clear that the sequence $(x_k)$ is
bounded, and hence it has a subsequence (without relabeling) that
converges weakly to $\ox\in\Omega$. Then (\ref{cover}) also holds true.

Now let $x\in \Omega$ and $t\geq 0$ satisfy $\Omega_i\subset x+tF$ for
$i=1,\ldots,n$.  Then $t\in \mathcal{E}$, and hence $r\leq t$. Thus
the smallest enclosing ball problem has a solution, and the proof is
now complete. $\h$

\begin{Proposition}\label{eq}  An element $\ox\in \Omega$ is a solution of the optimization problem (\ref{SEC}) with $r=C(\ox)$ if and only if $\ox$ is a solution of the smallest enclosing ball problem with smallest radius $r$.
\end{Proposition}
{\bf Proof: } Suppose that $\ox\in \Omega$ is a solution of the
optimization problem (\ref{SEC}) with $r=C(\ox)$. Then
\begin{equation*}
C_F(\ox; \Omega_i)\leq r \mbox{ for all } i=1,\ldots, n.
\end{equation*}
Thus $\Omega_i\subset \ox+rF$ for all $i=1,\ldots,n$. Let $x\in \Omega$ and
$t\geq 0$ satisfy
\begin{equation*}
\Omega_i\subset x+tF \mbox{ for all } i=1,\ldots,n.
\end{equation*}
Define $r^\prime= C(x)$. Then $r\leq r^\prime$. Moreover, it is also
clear that $r^\prime\leq t$. Thus $r\leq t$, and $\ox$ is a solution
of the smallest enclosing ball problem with smallest radius $r$.

Now suppose that $\ox\in \Omega$ is a solution of the smallest enclosing
ball problem with smallest radius $r$. We will prove that $C(\ox)=r$
and $C(\ox)\leq C(x)$ for all $x\in \Omega$ to show that $\ox$ is a
solution of problem (\ref{SEC}). Since $\Omega_i\subset \ox +rF$ for all
$i=1,\ldots,n$, one has $C(\ox)\leq r$. If $C(\ox)<r$, then choose
$r^\prime$ such that $C(\ox)<r^\prime <r$. Then $\Omega_i\subset \ox
+r^\prime F$ for all $i=1,\ldots,n$, which results in a
contradiction because $r$ is the smallest radius associated with
$\ox$. It follows that $r=C(\ox)$. For any $x\in \Omega$, define
$t=C(x)$. Then $\Omega_i\subset x+tF$ for all $i=1,\ldots,n$. Therefore,
$C(\ox)=r\leq t=C(x)$. The proof is now complete. $\h$
\vspace*{0.05in}

In what follows we will establish sufficient conditions for the
uniqueness of the smallest enclosing ball problem.

\begin{Theorem}\label{uniqueness1} Suppose that $X=\R^n$, $F$ is the Euclidean closed unit ball of $X$, and the
constraint $\Omega$ is a nonempty closed convex subset of $X$. Then
problem (\ref{SEC}) has a unique solution.
\end{Theorem}
{\bf Proof: }By Proposition \ref{eq}, in order to solve the smallest
enclosing ball problem, we only need to solve problem (\ref{SEC}).
The existence of an optimal solution under the assumptions made has
been proven in Theorem \ref{existence1}. Notice that $\ox$ is a
solution of the optimization problem (\ref{SEC}) if and only if it
is a solution of the following problem
%
\begin{align*}
\mbox{minimize } C^2(x) \mbox{ subject to } x\in \Omega
\end{align*}
where $C^2(x)=\max\{ [C_F(x;\Omega_i)]^2: i=1,\ldots,n\}$.

Since the maximum of a finite number of strictly convex functions on
$\Omega$ is a strictly convex function on this set, the proof reduces to
showing that each function
\begin{equation*}
c_i(x)=[C_F(x;\Omega_i)]^2=\sup_{\omega\in \Omega_i}||x-\omega||^2 \, \,
\mbox{ for } i = 1, \ldots, n
\end{equation*}
is strictly convex, where the definition of $C_F(x;\Omega_i)$ arises
from Proposition (\ref{p1}). It is obvious that the square norm
function is strictly convex on $X$. Fix $x, y\in \Omega$ with $x\neq y$
and $t\in (0,1)$. Denote $x_t=tx+(1-t)y \in \Omega$. Then there exists
$\omega_{t_i}\in \Omega_i$ for each $i$ such that
\begin{align*}
c_i(x_t)= ||x_t-\omega_{t_i}||^2=& ||tx+(1-t)y-\omega_{t_i}||^2\\
&=||t(x-\omega_{t_i})+(1-t)(y-\omega_{t_i})||^2\\
&< t||x-\omega_{t_i}||^2 +(1-t)||y-\omega_{t_i}||^2\\
&\leq tc_i(x)+(1-t)c_i(y).
\end{align*}
Therefore, each function $c_i(\cdot)$ is strictly convex on $\Omega$,
and thus $C^2(x)$ is strictly convex and a unique solution exists.
$\h$ \vspace*{0.05in}

The following examples show that the assumptions made in Theorem
\ref{uniqueness1} are essential.

\begin{Example}{\rm Let $X=\R^2$ and let $F$ be the Euclidean closed unit ball of $X$. Consider the smallest enclosing ball problem with the target set $\Omega_1=\{(0,0)\}$ and the constraint set \begin{equation*}
\Omega=\{(x_1,x_2): x_1^2+x_2^2=1\}.
\end{equation*}
The constraint set is nonconvex and in this case any point $x \in
\Omega$ is a solution of the smallest enclosing ball problem. }
\end{Example}

\begin{Example} {\rm Let $X=\R^2$ with unconstrained set $\Omega = \R^2$. Let $F=[-1,1]\times[-1,1]$. Define
\begin{equation*}
\Omega_1=\{ (0,1)\} \mbox{ and } \Omega_2=\{(0,-1)\}.
\end{equation*}
Then $\Omega_1$ and $\Omega_2$ are both convex. However, the unconstrained
smallest enclosing ball problem with target sets $\Omega_1$ and $\Omega_2$
has infinitely many solutions. In fact, any point of the set
\begin{equation*}
L=\{(x_1,x_2)\in X: x_1\in [-1,1], x_2=0\}
\end{equation*}
is a solution of the problem.}
\end{Example}

\section{The Smallest Intersecting Ball Problem}
\setcounter{equation}{0} Throughout this section we make the
following standing assumptions unless otherwise stated:

\emph{$X$ is a normed space; $F$ is a closed, bounded, convex set
that contains the origin as an interior point; the constraint set
$\Omega$ and the target sets $\Omega_i$, $i=1,\ldots,n$, are nonempty closed
subsets of $X$}.

\begin{Theorem}\label{existence2} Assume that one of the following statements holds: \\[1ex]
{\bf (i)} $X$ is finite dimensional, and one of the sets among $\Omega_i$, $i=1,\ldots,n,$ and $\Omega$ is bounded.  \\[1ex]
{\bf (ii)} $X$ is a reflexive Banach space, all of the sets $\Omega_i$,
$i=1,\ldots,n,$ and $\Omega$ are weakly closed, and at least one of them
is bounded.

Then the smallest intersecting ball problem has a solution. In this
case there exists $r\geq 0$ and $\ox\in \Omega$ such that
\begin{equation*}
(\ox +r F)\cap \Omega_i\neq \emptyset \mbox{ for all }i=1,\ldots, n,
\end{equation*}
and for any $x\in \Omega$ and $t\geq 0$ with $(x+tF)\cap \Omega_i\neq
\emptyset$ for all $i=1,\ldots,n$, one has $r\leq t$.
\end{Theorem}
{\bf Proof: }Consider the following set
\begin{equation*}
\mathcal{I}=\{ t\geq 0: \mbox{there exists } x\in \Omega \mbox{ with }
(x+tF)\cap \Omega_i\neq\emptyset \mbox{ for all }i=1,\ldots,n\}.
\end{equation*}
Then  $\mathcal{I}$ is nonempty since $0\in \mbox{int }F$. Moreover,
$\mathcal{I}$ is obviously bounded below. Let
\begin{equation*}
r=\inf \mathcal{I}\in [0,\infty).
\end{equation*}
Then there exists a sequence $(t_k)\subset \mathcal{I}$ that
converges to $r$. Let $(x_k)\subset \Omega$ satisfy
\begin{equation*}
(x_k +t_k F)\cap \Omega_i\neq \emptyset \mbox{ for all }i=1,\ldots,n.
\end{equation*}
Then there exist $f_{k,i}\in F$ and $\omega_{k, i}\in\Omega_i$, $k\in
\N$, $i=1,\ldots,n$, such that
\begin{equation*}
x_k+t_kf_{k,i}=\omega_{k,i} \mbox{ for all } k \mbox{ for all
}i=1,\ldots, n.
\end{equation*}
Let us focus on case {\bf (i)}. We will first show that $(x_k)$ is
bounded under the assumptions made. Without loss of generality,
suppose that $\Omega_1$ is bounded. One has
\begin{equation*}
x_k+t_kf_{k,1}=\omega_{k,1} \mbox{ for all } k.
\end{equation*}
Thus
\begin{equation*}
||x_k||\leq t_k||f_{k,1}||+||\omega_{k,1}||.
\end{equation*}
Since both $F$ and $\Omega_1$ are bounded, $(x_k)$ is a bounded
sequence; thus there exists a convergent subsequence $(x_k)$
(without relabeling). Since $F$ is closed and bounded, and $\Omega$ is closed, we can assume that $f_{k,i}\to
f_i\in F$ for $i=1,\ldots,n$ and $x_k \to \ox\in \Omega$. For
$i=1,\ldots,n$, one also has
\begin{equation*}
x_k +t_k f_{k,i} =\omega_{k,i}\to \ox+r f_i \mbox{ as } k\to \infty,
\end{equation*}
and $\ox+r f_i \in \Omega_i$ since each $\Omega_i$ is a closed set.
Moreover, $\ox+r f_i \in \ox + rF$ for each $i=1, \ldots, n$. Thus
$(\ox+ r F)\cap \Omega_i\neq \emptyset$. For any $x\in \Omega$ and $t\geq 0$
with $(x+tF)\cap \Omega_i\neq \emptyset$ for all $i=1,\ldots,n$, one has
$t\in \mathcal{I}$. Thus $r\leq t$.

The proof of the result under case {\bf (ii)} is similar where the
weak convergence of $(x_{k}), (f_{k,i})$ and $(w_{k,i})$ for $i = 1,
\ldots, n$ are  taken into account. The proof is now complete. $\h$

\begin{Lemma}\label{minimum}  Suppose one of the following:\\[1ex]
{\bf (i)} $X$ is finite dimensional, and $Q$ is a closed subset of $X$;\\[1ex]
{\bf (ii)} $X$ is reflexive, and $Q$ is a nonempty, weakly, closed subset of $X$. \\[1ex]
Let
\begin{equation*}
\bar t =T_F(\ox; Q)=\inf \{ t\geq 0: (\ox+t F)\cap  Q \neq \emptyset \}.
\end{equation*}
Then $\bar t\geq 0$ and $(\ox +\bar t F)\cap Q \neq \emptyset$.
\end{Lemma}
{\bf Proof: } Assume case (i) where  $X$ is finite dimensional, and
let $\bar x \in X$. It is clear that $T_F(\ox; Q)$ is a finite
number. Let $(t_k)$ be a sequence of nonnegative integers that
converges to $\bar t \geq 0$, where
\begin{equation*}
(\ox +t_k F)\cap Q\neq \emptyset \mbox{ for all }k.
\end{equation*}
Let $f_k\in F$ and $q_k\in Q$ satisfy that $\ox+t_k f_k=q_k$. Since
$F$ is closed and bounded,  we can assume without loss of generality
that $f_k\to f\in F$. Then
\begin{equation*}
\ox+t_k f_k=q_k\to \ox+\bar tf \in Q.
\end{equation*}
Therefore, $(\ox+\bar t F)\cap Q\neq \emptyset.$

Similarly, we can prove case (ii) where $X$ is reflexive and $Q$ is
a nonempty, weakly, closed subset of $X$. $\h$

\begin{Proposition} Suppose one of the following: \\[1ex]
{\bf (i)} $X$ is finite dimensional;\\[1ex]
{\bf (ii)} $X$ is reflexive and $\Omega$ and $\Omega_i$, $i=1,\ldots,n$, are nonempty, weakly, closed subsets of $X$.\\[1ex]
Then  $\ox\in \Omega$ is a solution of the optimization problem
(\ref{SIB}) with $r=T(\ox)$ if and only if $\ox\in \Omega$ is a solution
of the smallest intersecting ball problem with smallest radius $r$.
\end{Proposition}
{\bf Proof: } Suppose that $\ox\in\Omega$ is a solution of the
optimization problem (\ref{SIB}) and $r=T(\ox)$. Let
\begin{equation*}
t_i=T_F(\ox; \Omega_i) \mbox{ for } i=1,\ldots, n.
\end{equation*}
This implies $t_i \leq r$ for all $i=1,\ldots, n$. By
Lemma (\ref{minimum}),
\begin{equation*}
(\ox+t_i F)\cap \Omega_i\neq \emptyset \mbox{ for all }i=1,\ldots,n.
\end{equation*}
Thus $(\ox +rF)\cap \Omega_i\neq \emptyset$, which follows from the fact
that $t_iF\subset rF$ for all $i=1,\ldots,n$ under the assumptions that
$F$ is convex and $0\in F$.

Now suppose that $x\in \Omega$ and $t\geq 0$ satisfy
\begin{equation*}
(x+tF)\cap \Omega_i\neq \emptyset \mbox{ for all }i=1,\ldots,n.
\end{equation*}
Then $T_F(x;\Omega_i)\leq t$ for all $i=1,\ldots, n$. This implies
$r=T(\ox)\leq T(x)\leq t.$ Thus $\ox$ is a solution of the smallest
intersecting ball problem with radius $r$.

Conversely, suppose that $\ox\in \Omega$ is a solution of the smallest
intersecting ball problem with  smallest radius $r\geq 0$. We will
prove that $\ox$ is a solution of the optimization problem
(\ref{SIB}) and $r=T(\ox)$.  One has that
\begin{equation*}
(\ox+rF)\cap \Omega_i\neq \emptyset \mbox{ for all }i=1,\ldots,n.
\end{equation*}
This implies
\begin{equation*}
T(\ox)=\max\{ T_F(\ox;\Omega_i): i=1,\ldots,n\}\leq r.
\end{equation*}
If $T(\ox)<r$, then there exists a real number $s$ such that $T(\ox)<s <r$ and
we easily see that $(\ox+sF)\cap \Omega_i\neq \emptyset$ for
$i=1,\ldots,n$, which contradicts the minimal property of $r$. Thus
$T(\ox)=r$. Now take any $x\in \Omega$. Define $r'=T(x)$. Then
$(x+r'F)\cap \Omega_i\neq\emptyset$ for all $i=1,\ldots,n$. Thus $r\leq
r'$ or equivalently $T(\ox)\leq T(x)$. Therefore, $\ox$ is a
solution of (\ref{SIB}). The proof is complete. $\h$\vspace*{0.05in}

The following theorem provides natural sufficiency conditions
guaranteeing the uniqueness of the solution for the smallest
intersecting ball problem.

\begin{Theorem}\label{uniqueness2} Let $X$ be a Hilbert space and let $F$ be the closed, unit ball of $X$.
Suppose that $\Omega$ is a nonempty, closed, convex set, $\Omega_i$, $i=1,
\ldots, n$, are strictly convex, and at least one of the sets among
$\Omega_i$, $i=1,\ldots, n$ and $\Omega$ is bounded. Suppose further that
\begin{equation}\label{emptyassump}
\cap_{i=1}^n[\Omega_i\cap\Omega]=\emptyset.
\end{equation}
Then the optimization problem (\ref{SIB}) has a unique solution.
\end{Theorem}
{\bf Proof: }The existence of an optimal solution follows from
Theorem \ref{existence2}.

Since $F$ is the closed unit ball of $X$, the minimal time function
$T_F(\cdot;\Omega_i)$ reduces to the distance function $d(\cdot; \Omega_i)$
for $i=1,\ldots,n$. Moreover, $\ox$ is a solution of the smallest
intersecting ball problem if and only if it is a solution to the
optimization problem (\ref{SIB}) by Proposition 3.3.

For $x \in \Omega$, consider the function
\begin{equation*}
S(x)=\mbox{max}\{ [d(x;\Omega_i)]^2: i=1, \ldots, n\}.
\end{equation*}
Then $\ox$ is a solution of problem (\ref{SIB}) if and only if $\ox$
is a solution of the problem
\begin{equation*}
\mbox{ minimize } S(x) \mbox{ subject to } x\in \Omega.
\end{equation*}
We will prove that $S$ is strictly convex on $\Omega$. Fix $x, y\in \Omega$
with $x\neq y$ and $t\in (0,1)$. Denote $x_t=tx+(1-t)y \in \Omega$. Let
$z \in \Omega$ and define $I(z):=\{i=1,\ldots, n: T(z)=T_F(z; \Omega_i)\}$.
Then for $i \in I(x_t)$ we have $T(x_t)= d(x_t;\Omega_i)$. Let $u, v\in
\Omega_i$  satisfy that $d(x;\Omega_i)=||x-u||$ and $d(y;\Omega_i)=||y-v||$. It
follows that
\begin{align*}
S(x_t)&= [d(x_t;\Omega_i)]^2\\
&\leq ||x_t-(tu+(1-t)v)||^2\\
&=||t(x-u)+ (1-t)(y-v)||^2 \, \, \mbox{ by definition of } x_t,\\
& = t^2 ||x-u||^2 +2t(1-t) \la x-u, y-v\ra + (1-t)^2 ||y-v||^2\\
&\leq t^2 ||x-u||^2 + 2t(1-t) ||x-u||\cdot ||y-v|| +(1-t)^2 ||y-v||^2\\
&\leq t^2 ||x-u||^2 + t(1-t) (||x-u||^2 + ||y-v||^2) +(1-t)^2 ||y-v||^2\\
&= t ||x-u||^2 + (1-t) ||y-v||^2\\
&=t [d(x;\Omega_i)]^2 + (1-t) [d(y;\Omega_i)]^2\\
&\leq t S(x)+ (1-t) S(y).
\end{align*}
Thus $S(x)$ is convex on $\Omega$. Moreover,
\begin{align*}
S(x_t)&= [d(x_t;\Omega_i)]^2 \leq  [td(x;\Omega_i)+ (1-t) d(y;\Omega_i)]^2\\
&\leq t^2 ||x-u||^2 + t(1-t) (||x-u||^2 + ||y-v||^2) +(1-t)^2 ||y-v||^2\\
&= t ||x-u||^2 + (1-t) ||y-v||^2\\
&=t [d(x;\Omega_i)]^2 + (1-t) [d(y;\Omega_i)]^2\\
&\leq t S(x)+ (1-t) S(y).
\end{align*}
Now, suppose the equality $S(x_t)= t S(x)+ (1-t) S(y)$ holds; we
will show this leads to a contradiction. We have
\begin{equation*}
\la x-u,y-v\ra=||x-u||.||y-v||, ||x-u||=||y-v|| \mbox{ and }
d(tx+(1-t)y; \Omega_i)=td(x;\Omega_i)+(1-t) d(x;\Omega_i),
\end{equation*}
which imply that
\begin{equation*}
x-u=y-v \mbox{ and } d(tx+(1-t)y; \Omega_i)=td(x;\Omega_i)+(1-t) d(x;\Omega_i).
\end{equation*}
This implies $u\neq v$ and $d(tx+(1-t)y;\Omega_i)=d(x;\Omega_i)=d(y;\Omega_i)$.
Since $T(x_t)=d(x_t; \Omega_i)>0$ by (\ref{emptyassump}), one has $x_t\notin \Omega_i$. Using the
strict convexity of $\Omega_i$, one has $tu+(1-t)v\in \mbox{ int }\Omega_i$.
Thus
\begin{equation*}
d(x_t;\Omega_i)< ||x_t - [tu+(1-t)v]|| \leq t||x-u||+(1-t)||y-v||=
d(x;\Omega_i).
\end{equation*}
Indeed, let $\delta>0$ satisfy $\B(tu+(1-t)v; \delta)\subset \Omega_i$.
Denote $c=tu+(1-t)v$. Then $c+ \delta \dfrac{x_t-c}{||x_t-c||}\in
\Omega_i$ and
\begin{equation*}
||x_t- \delta \dfrac{x_t-c}{||x_t-c||}|| = ||x_t-c||-\delta
<||x_t-c||,
\end{equation*}
which is a contradiction. Therefore, $S$ must be strictly convex, and the problem has a unique solution. The proof is now complete. $\h$

The following two examples illustrate the need for the assumptions
in the above theorem that provides the uniqueness of the solution to
the smallest intersecting ball problem.

\begin{Example} {\rm Let $X=\R^2$ with $\Omega = X$ and let $F$ be the Euclidean closed unit ball in $X$. Define target sets
\begin{equation*}
\Omega_1=\{(x_1,x_2)\in X: x_2\geq 1\} \mbox{ and } \Omega_2=\{(x_1,x_2)\in
X: x_2\leq -1\}.
\end{equation*}
 Note that $\Omega_1$ and $\Omega_2$ are both convex and violate the strict convexity assumptions in Theorem \ref{uniqueness2}. In addition, none of the target sets $\Omega_1$ and $\Omega_2$ or $\Omega$ is bounded.  In this case, the unconstrained smallest intersecting ball problem  has infinitely many solutions. In fact, any point of the set
\begin{equation*}
K=\{(x_1,x_2)\in X: x_2=0\}
\end{equation*}
is a solution of the problem.}
\end{Example}

\begin{Example} {\rm Let $X=\R^2$ with $\Omega = X$. Let $F=[-1,1]\times[-1,1]$ be the unit square centered at the origin; thus $F$ violates the assumptions of Theorem \ref{uniqueness2}. Define
\begin{equation*}
\Omega_1=\B((0,2); 1) \mbox{ and } \Omega_2=\B((0,-2);1).
\end{equation*}
Then $\Omega_1$ and $\Omega_2$ are both strictly convex and bounded.
However, the unconstrained smallest intersecting ball problem with
target sets $\Omega_1$ and $\Omega_2$ has infinitely many solutions. In
fact, any point of the set
\begin{equation*}
L=\{(x_1,x_2)\in X: x_1\in [-1,1], x_2=0\}
\end{equation*}
is a solution of the problem.}
\end{Example}

\section{Conclusions}

This paper is a part of our project involving \emph{set facility
location problems}. The main idea is to consider a much broader
situation where singletons in the classical models of facility
location problems are replaced by sets. The new extension seems to
be interesting for both the theory and applications to various
location models, optimal networks, wireless communications, etc.
Moreover, it sheds new lights on classical geometry problems.

Our next goal is to study optimality conditions and numerical
algorithms for the smallest enclosing ball problem and the smallest
intersecting ball problem. Based on the approach we have developed
in \cite{mnft,naj,mnj1,mnj2}, we foresee the potential of success of
this future work.


\small


\begin{thebibliography}{99}


\bibitem{chm} Cheng, D., Hu, X., Martin, C.: On the smallest enclosing balls, Commun. Inf. Syst. 6 (2006), 137--160.

\bibitem{DZ} Deviille, R., Zizler, V.E.: Farthest points in $w^*-$ compact sets, Bull. Austral. Math. Soc. 38 (1988), 433--439.


\bibitem{ljc} Drager L., Lee. J., Martin, C.: On the geometry of the smallest circle enclosing a finite set of points,
J. Franklin Inst. 344 (2007), 929--940.


\bibitem{Lau} Lau, K.S.: Farthest points in weakly compact sets, Israel J. Math. 22 (1975), 168--174.
I. Fundamentals, Springer-Verlag, Berlin (1993).

\bibitem{mor06a} Mordukhovich, B.S.: Variational Analysis and
Generalized Differentiation, I: Basic Theory, II: Applications,
Grundlehren Series (Fundamental Principles of Mathematical
Sciences), Vols. 330 and 331, Springer, Berlin (2006).


\bibitem{bmn10} Mordukhovich, B.S., Nam, N.M.: Subgradients of minimal time functions
under minimal assumptions. J. Convex Anal. 18 (2011), 915--947.

\bibitem{mnft} Mordukhovich, B.S., Nam, N.M.: Applications of variational analysis to
a generalized Fermat-Torricelli problem. J. Optim. Theory Appl. 148
(2011), 431--454.

\bibitem{mnj1} Mordukhovich, B.S., Nam, N.M., Salinas, J.: Applications of variational analysis to a generalized Heron
problemroblem, to appear in Applicable Analysis.

\bibitem{mnj2} Mordukhovich, B.S., Nam, N.M., Salinas, J.: Applications of variational analysis to a generalized Heron
problemroblem, to appear in American Math Monthly.

\bibitem{naj} Nam, N.M., An N.T., Salinas, J.: Applications of convex analysis to
the smallest intersecting ball problem, to appear in Journal of
Convex Analysis.

\bibitem{frank} Nielsen, F., Nock, R.: Approximating smallest enclosing balls with applications to machine
learning.  Internat. J. Comput. Geom. Appl. 19 (2009), 389--414.


\bibitem{syl} Sylvester, J.J.: A question in the geometry of situation. Quarterly Journal of Pure and Applied Mathematics 1:79 (1857).


\bibitem{wel} Welzl, E.: Smallest enclosing disks (balls ellipsoids).
H. Maurer, editor, Lecture Notes in Comput. Sci. 555 (1991),
359--370.

\bibitem{WS} Westphal, U., Schwartz, T.: Farthest points and monotone operators, Bull. Austral. Math. Soc. 58 (1998), 75--92.






\end{thebibliography}
\end{document}